# Kemeny's Function for Markov Chains and Markov Renewal Processes

Jeffrey J. Hunter

*Department of Mathematical Sciences*
*School of Engineering, Computer and Mathematical Sciences,*
*Auckland University of Technology, Auckland, New Zealand.*
*Email: jeffrey.hunter@aut.ac.nz*

February 23, 2018

## Abstract

Extensions of Kemeny's constant, as derived for irreducible finite Markov chains in discrete time, to Markov renewal processes and Markov chains in continuous time are discussed. Three alternative Kemeny's functions and their variants are considered. Typically, they lead to a constant if and only if the mean holding times between the states in the Markov renewal process are constant. However one particular variant leads to a constant, analogous to the discrete time Markov chain result. Specifically, if the state space is finite, the weighted sum of the mean first passage times (omitting the mean return time) with the stationary probabilities associated with the continuous time semi-Markov process is a constant for any Markov renewal process. Expressions for the Kemeny's functions and the relevant constants are derived for Markov renewal processes and special cases involving continuous time Markov chains and birth and death processes.

## 1. Introduction

A key property of irreducible finite discrete time Markov chains (*MC's*) is the existence of Kemeny's constant. The result is simple in form. Let the *MC* have state space *S,* stationary probabilities $\{\pi_i\}$ $(i \in S)$ and mean first passage times $m_{ij}$ $(i,j \in S)$ (mean recurrence time when $i = j$). The expression $\sum_{j \in S} \pi_j m_{ij}$, Kemeny's function, $k_i$, is in fact a constant, Kemeny's constant, not depending on $i$. This constant was identified and explored by John Kemeny and first appears in Kemeny and Snell (1960). It has the simple interpretation as the hitting time $T$ to a 'random" state chosen from the stationary distribution starting from any fixed state $i$ but has the property that it is independent of the particular state used as the starting point. For various derivations and applications of this constant to different situations see Hunter (2014). Catral et al (2010) also identified that Kemeny's constant is the same as the expected time to mixing, i.e., time to reach stationarity in an ergodic *MC*, which was also explored in Hunter (2006)).

One major interest has been efforts to give a physical interpretation, as opposed to mathematical arguments, as to why this expression should be constant. Hunter (2014) chronicles efforts to find a simple physical interpretation that were initially stimulated by a prize given by John Kemeny, as described in Snell (1975). The first recorded attempt appears to have been given by Peter Doyle in 1983, using the maximum principle, as highlighted in Grinstead and Snell (1997), when the question is posed as to whether Peter should have been given the prize. Further recent arguments are given by Gustafson and Hunter (2016), and Bini, Hunter, Meini, Latouche and Taylor (2017).

What happens if we broaden the context and consider the equivalent expression for Markov renewal processes (*MRP's*)? Does Kemeny's function lead to a constant in this environment? If not, under what conditions are required for it to be constant? What about continuous time *MC's*?

We first summarize a derivation of Kemeny's constant for finite discrete time *MC's*, followed by derivations of related Kemeny's functions for *MRP's*. Some variants analogous to omitting the $\pi_j\, m_{jj}$ term are also explored. The conditions for constancy and the extension to *MC's* in continuous time are then explored.

## 2. Kemeny's function for discrete time Markov chains on a finite state space.

Let $\{X_n\}$ be a finite irreducible discrete time *MC* with state space $S = \{1, 2, \ldots, m\}$ and transition matrix $P = [p_{ij}]$, where $p_{ij} = P[X_{n+1} = j \mid X_n = i]$. Then a stationary distribution $\{\pi_i\}$ exists for $X_n$, where $\pi_j = \sum_{i=1}^{m} \pi_i p_{ij}$ for all $j \in S$ with $\sum_{i=1}^{m} \pi_i = 1$, (Feller (1950)).

Let $\boldsymbol{\pi}^T = (\pi_1, \pi_2, \ldots, \pi_m)$ be the stationary probability vector of $\{X_n\}$ and $\boldsymbol{\pi}^T$ is the unique left-eigenvector corresponding to the eigenvalue 1, i.e.

$$\boldsymbol{\pi}^T (I - P) = \boldsymbol{0}^T \text{ and } \boldsymbol{\pi}^T \boldsymbol{e} = 1, \tag{2.1}$$

where $\boldsymbol{e}^T = (1,1,\ldots,1)$ is a vector of 1's of dimension $m$.

It is also the case that for the eigenvalue 1, $P$ has a right-eigenvector that is unique up to multiplicative constant of the vector $\boldsymbol{e}$. i.e.

$$(I - P)\boldsymbol{x} = \boldsymbol{0} \Leftrightarrow \boldsymbol{x} = c\boldsymbol{e}, \tag{2.2}$$

for some real $c$.

Let $N_{ij}$ denote the generic first-passage time random variable for the *MC* $X_n$ to pass from state $i$ to state $j$, (return time when $i = j$), so that $N_{ij} = \inf\{n \geq 1 : X_n = j \mid X_0 = i\}$. Under the assumption of irreducibility, the $N_{ij}$ are proper random variables and the $m_{ij} = \mathrm{E}\,(N_{ij} \mid X_0 = i)$ are well defined and finite for all $i, j \in S$.

It is well known, by a simple first-passage decomposition argument, (e.g. Kemeny and Snell (1960)) that for all $i, j \in S$,

$$m_{ij} = 1 + \sum_{k \neq j} p_{ik} m_{kj}. \tag{2.3}$$

Using the notation $M = [m_{ij}]$, $M_d = [\delta_{ij} m_{ij}] = diag(m_{11}, \ldots, m_{mm})$, $E = [1] = \boldsymbol{e}\boldsymbol{e}^T$, equations (2.3) can be re-expressed as

$$(I - P)M = E - PM_d. \tag{2.4}$$

In particular, it is a well known standard result for *MC's* (Kemeny and Snell (1960)), that

$$m_{ii} = 1/\pi_i, \tag{2.5}$$

so that if $\Pi = \boldsymbol{e}\boldsymbol{\pi}^T$ then $D = M_d = (\Pi_d)^{-1}$.

Hunter (1982), (see also Hunter (1983)), gave a general solution to equation (2.4), using any one-condition generalized inverse $G$ of $I - P$, i.e. any $G$ such that $(I - P)G(I - P) = I - P$, that

$$M = [G\Pi - E(G\Pi)_d + I - G + EG_d]D. \tag{2.6}$$

We define Kemeny's function as $k_i \equiv \sum_{j=1}^{m} m_{ij}\pi_j$ and thus if $\boldsymbol{k}^T = (k_1, k_2, ..., k_m)$ then $\boldsymbol{k} = M\boldsymbol{\pi}$.
Since, from (2.4), $M_d\boldsymbol{\pi} = \boldsymbol{e}$, we have, from (2.5),

$$(I - P)\boldsymbol{k} = (I - P)M\boldsymbol{\pi} = \boldsymbol{e}\boldsymbol{e}^T\boldsymbol{\pi} - P\boldsymbol{e} = \boldsymbol{e} - \boldsymbol{e} = \boldsymbol{0}. \tag{2.7}$$

Thus, from (2.2), it follows that $\boldsymbol{k} = k\boldsymbol{e}$ for some constant $k$, implying that $k_i = k$ for all $i$.
So that Kemeny's function is in fact a constant for all $i$, known as Kemeny's constant $K_C$.

From (2.6) (or Hunter (2006)) it is easily shown that if $G$ is any g-inverse of $I - P$, then, since $\boldsymbol{e}^T G_d \boldsymbol{e} = tr(G)$ and $\boldsymbol{e}^T(G\Pi)_d\boldsymbol{e} = tr(G\Pi)$,

$$K_C = 1 + tr(G) - tr(G\Pi). \tag{2.8}$$

In elemental form, if $G = [g_{ij}]$ with $g_{j.} = \sum_{k=1}^{m} g_{jk}$,

$$K_C = 1 + \sum_{j=1}^{m}(g_{jj} - g_{j.}\pi_j). \tag{2.9}$$

In particular, $$K_C = tr(Z) = 1 + tr(A^{\#}), \tag{2.10}$$

where $Z = [I - P + \Pi]^{-1}$ is the fundamental matrix of $I - P$ and $A^{\#} = Z - \Pi$ is the group inverse of $I - P$.

An alternative expression for Kemeny's constant can be given in terms of eigenvalues of $P$, $\{\lambda_i, i = 1, ..., m\}$ which, since $P$ is irreducible, are such that $\lambda_1 = 1$, with $|\lambda_i| \le 1$ and $\lambda_i \ne 1$, ($i$ = 2, .., $m$) (see Hunter (2006), Hunter (2014)):

$$K_C = 1 + \sum_{j=2}^{m} 1/(1 - \lambda_j). \tag{2.11}$$

See also Catral et al (2010) for other alternative expressions of $K_C$.

An alternative variant of Kemeny's constant is given by $K^\circ_C \equiv \sum_{j=1, j\ne i}^{m} \pi_j m_{ij}$, since, from (2.5), $m_{ii}\pi_i = 1$, so that $K^\circ_C = K_C - 1$.

The key observation is that for a discrete time Markov chain, the function $k_i \equiv \sum_{j=1}^{m} \pi_j m_{ij}$, or the variant $k^\circ_i \equiv \sum_{j=1, j\ne i}^{m} \pi_j m_{ij}$, are each a constant, independent of the initial state $i$, given by $K_C$ and $K^\circ_C$, respectively, and can be found in a variety of ways.

## 3. Kemeny's function for Markov renewal processes on a finite state space.

Let us set the scene with a review of the definitions used in Markov renewal theory.
Let $\{(X_n, T_n), n \geq 0\}$ be a finite irreducible *MRP* with underlying jump *MC* $\{X_n, n \geq 0\}$ having transition matrix $P$ and state space $S = \{1, 2, \ldots, m\}$. The semi-Markov kernel is given by $Q_{ij}(t) = P\{X_{n+1} = j, T_{n+1} - T_n \leq t | X_n = i\}$ so that $P = [p_{ij}] = [Q_{ij}(+\infty)]$.
Let $F_{ij}(t)$ be the distribution function of the length of time between transitions of the *MRP*, given that the process makes a transition from state $i$ to state $j$, i.e.
$F_{ij}(t) = P\{T_{n+1} - T_n \leq t | X_n = i, X_{n+1} = j\}$ so that, if $p_{ij} > 0$, $Q_{ij}(t) = p_{ij}F_{ij}(t)$.
Let $\mu_{ij} = \int_0^\infty t\,dQ_{ij}(t) = p_{ij}\int_0^\infty t\,dF_{ij}(t) = p_{ij}E[T_{n+1} - T_n \mid X_n = i, X_{n+1} = j]$ where $E[T_{n+1} - T_n \mid X_n = i, X_{n+1} = j]$ is the mean holding time in state $i$ prior to moving to state $j$, given that a transition takes place between states $i$ and $j$. We assume that these first moments are all finite. Let $P^{(1)} = \left[\mu_{ij}\right]$ and let $\mu_i$ be the mean sojourn time of the *MRP* on any visit to state $i$, i.e. $\mu_i = \sum_k \mu_{ik}$ so that $P^{(1)}\boldsymbol{e} = \boldsymbol{\mu}$, where $\mu^T = (\mu_1, \mu_2, \ldots, \mu_m)$.

Let $T_{ij}$ be the time for a first passage from state $i$ to state $j$ in the *MRP*.
Let $G_{ij}(t)$ be the distribution function of $T_{ij}$ and let $m_{ij} = \int_0^\infty t\,dG_{ij}(t) = E[T_{ij} \mid X_0 = i]$, be the mean first passage time from state $i$ to state $j$ in the *MRP*.

By a first passage decomposition, (See Pyke (1961), Hunter (1969)),

$$G_{ij}(t) = Q_{ij}(t) + \sum_{k \neq j}\int_0^\infty Q_{ik}(t-u)\,dG_{kj}(u), \tag{3.1}$$

which leads (Corollary 2.1.1, Hunter (1969)) to relationships between the first moments

$$m_{ij} = \mu_i + \sum_{k \neq j} p_{ik}m_{kj}. \tag{3.2}$$

There are two stationary distributions associated with the *MRP* $\{(X_n, T_n), n \geq 0\}$ – one $\{\pi_i\}$ associated with the discrete time *MC* $\{X_n, n \geq 0\}$ and one $\{\varpi_j\}$ associated with the minimal semi-Markov process, *SMP*, $\{X(t), t \geq 0\}$ where $X(t) = X_n$ for $T_n \leq t < T_{n+1}$. (i.e. $X(t)$ is the state that the *MRP* is occupying at time *t*.) If the *MC* is aperiodic $\pi_j = \lim_{n\to\infty} P\{X_n = j | X_0 = i\}$, while $\varpi_j = \lim_{t\to\infty} P\{X(t) = j | X_0 = i\}$ or, if the *SMP* is stationary, $\varpi_j = \lim_{t\to\infty} P\{X(t) = j\}$.

Let $\boldsymbol{\pi}^T = \{\pi_1, \pi_2, \ldots, \pi_m\}$ and $\boldsymbol{\varpi}^T = \{\varpi_1, \varpi_2, \ldots, \varpi_m\}$ be the stationary probability vectors associated with the *MC* $\{X_n\}$ and the *SMP* $\{X(t)\}$, respectively.

Both of these sets of probabilities or vectors are interrelated. Let $\lambda \equiv \boldsymbol{\pi}^T\boldsymbol{\mu}$ and $\Lambda \equiv diag(\mu_1, \mu_2, \ldots, \mu_m)$. ($\lambda$ has been called the "mean asymptotic increment", (Keilson and Wishart (1964). See also Hunter (1969)).
If $\{X_n\}$ is irreducible and aperiodic, $\mu_i < \infty$ for all $i \in S$, then (Çinlar (1975)),

$$\boldsymbol{\varpi}^T = \boldsymbol{\pi}^T \Lambda/\lambda \text{ and } \boldsymbol{\pi}^T = \lambda\boldsymbol{\varpi}^T\Lambda^{-1}, \tag{3.3}$$

or, equivalently, $\boldsymbol{\varpi} = \Lambda\boldsymbol{\pi} / \lambda$ and $\boldsymbol{\pi} = \lambda\Lambda^{-1}\boldsymbol{\varpi}$.

We have various possible ways of specifying Kemeny's function for *MRP's*. The natural approach is to specify Kemeny's function, analogous to the definition for *MC's*. We consider the alternative variants later.

**Definition 1**: Let $k_i^{(1)} \equiv \sum_{j=1}^{m} \pi_j m_{ij}$, where $\pi_j$ is the stationary probability of being in state $j$ in the embedded discrete-time *MC* $\{X_n\}$, be a Kemeny's function for a finite irreducible $m$-state *MRP* with mean first passage times $m_{ij}$, (mean recurrence time when $i = j$).

The rationale for this definition is reflected in the observation that Kemeny's constant in a *MC* is equivalent to the expected time to mixing – i.e. starting at any state $i$ the time to reach a state selected from the stationary distribution of the *MC*. In a *MRP* this time must occur when the *MRP* moves to a particular state for the first time and that is governed by the underlying embedded *MC*.

**Definition 2**: Let $k_i^{(2)} \equiv \sum_{j=1}^{m} \varpi_j m_{ij}$, where $\varpi_j$ is the stationary probability of being in state $j$ in the continuous time *SMP* $\{X(t)\}$, be a Kemeny's function for a finite irreducible $m$-state *MRP* with mean first passage times $m_{ij}$, (mean recurrence time when $i = j$)

Thus we use the limiting (in the case of aperiodic *MRP*'s) or the continuous time stationary state distribution rather than that of the embedded *MC* to select the target state.

Definition 2 was used in effect in Bini et. al. (2017) in the setting of Markov chains in continuous time (their Section 4).

In the case of discrete time MC's observe, by virtue of (2.5) that $K_C = \sum_{j=1}^{m} m_{ij}/m_{jj}$. This leads to another possible definition of Kemeny's function for an *MRP*.

**Definition 3**: Let $k_i^{(3)} \equiv \sum_{j=1}^{m} m_{ij}/m_{jj}$, where $m_{ij}$ is the mean first passage time between states $i$ and $j$ (mean recurrence time when $i = j$) in the *MRP* $\{X_n, T_n\}$, be a Kemeny's function for a finite irreducible $m$-state *MRP* with mean first passage times $m_{ij}$.

Let $\boldsymbol{k}^{(l)T} = (k_1^{(l)}, k_2^{(l)}, ..., k_m^{(l)})$ for $l$ = 1, 2, and 3. The vector versions of the respective Kemeny functions yield $\boldsymbol{k}^{(1)} = M\boldsymbol{\pi}$, $\boldsymbol{k}^{(2)} = M\boldsymbol{\varpi}$, and $\boldsymbol{k}^{(3)} = M(M_d)^{-1}\boldsymbol{e}$.

The first two definitions are interrelated through (3.3) with $\boldsymbol{k}^{(1)} = \lambda M\Lambda^{-1}\boldsymbol{\varpi}$ and $\boldsymbol{k}^{(2)} = M\Lambda\boldsymbol{\pi} / \lambda$.

Let $M = [m_{ij}]$ be the mean first passage time matrix of a finite irreducible *MRP* $\{(X_n,T_n), n \geq 0\}$. Let $M_d = [\delta_{ij}m_{ij}]$, with $\boldsymbol{e}^T$ and $E$ as defined in Section 2. The following result follows from (3.2) and appears in Corollary 2.3.1, Hunter (1969) and also in Section 5.2, Hunter (1982): $M$ satisfies the matrix equation

$$(I - P)M = P^{(1)}E - PM_d. \tag{3.4}$$

Pre-multiplication of (3.4) by $\boldsymbol{\pi}^T$ implies $\boldsymbol{\pi}^T P^{(1)} E = \boldsymbol{\pi}^T M_d$. Since $\boldsymbol{\mu} = P^{(1)}\boldsymbol{e}$ and $E = \boldsymbol{ee}^T$ we thus deduce that $\boldsymbol{\pi}^T \boldsymbol{\mu e}^T = \boldsymbol{\pi}^T M_d$. Since $\lambda \equiv \boldsymbol{\pi}^T \boldsymbol{\mu}$ we have that $\lambda \boldsymbol{e}^T = (\pi_1 m_{11},....,\pi_m m_{mm})$. Thus for all $i$,

$$\lambda = \pi_i m_{ii}. \tag{3.5}$$

Further $D \equiv M_d = \lambda(\Pi_d)^{-1}$ where $\Pi = \boldsymbol{e\pi}^T$, and consequently

$$D\boldsymbol{\pi} = \lambda diag(1/\pi_1,...,1/\pi_m)\boldsymbol{\pi} = \lambda \boldsymbol{e}. \tag{3.6}$$

**Theorem 1:** For a finite irreducible *MRP* with underlying irreducible jump *MC* having transition matrix $P$ with stationary probability vector $\boldsymbol{\pi}$, mean sojourn time vector $\boldsymbol{\mu}$ and $\lambda \equiv \boldsymbol{\pi}^T \boldsymbol{\mu}$, Kemeny's function vector, $\boldsymbol{k}^{(1)}$ under Definition 1, can be expressed as follows:

(*a*) Let $G$ be any g-inverse of $I - P$, then

$$\boldsymbol{k}^{(1)} = G\boldsymbol{\mu} - tr(G\boldsymbol{\mu\pi}^T)\boldsymbol{e} + \lambda\boldsymbol{e} - \lambda G\boldsymbol{e} + \lambda tr(G)\boldsymbol{e}. \tag{3.7}$$

(*b*) If $Z = [I - P + \Pi]^{-1}$ is the fundamental matrix of the embedded *MC*,

$$\boldsymbol{k}^{(1)} = Z\boldsymbol{\mu} - tr(Z\boldsymbol{\mu\pi}^T)\boldsymbol{e} + \lambda tr(Z)\boldsymbol{e}. \tag{3.8}$$

(*c*) If $A^{\#}$ is the group inverse of $I - P$,

$$\boldsymbol{k}^{(1)} = A^{\#}\boldsymbol{\mu} - tr(A^{\#}\boldsymbol{\mu\pi}^T)\boldsymbol{e} + \lambda\boldsymbol{e} + \lambda tr(A^{\#})\boldsymbol{e}. \tag{3.9}$$

(*d*) If $\widetilde{G} \equiv [I - P + \boldsymbol{\mu u}^T]^{-1}$ (with $\boldsymbol{u}^T\boldsymbol{e} \neq 0$)

$$\boldsymbol{k}^{(1)} = \lambda\left[I - \widetilde{G} + tr(\widetilde{G})\right]\boldsymbol{e}. \tag{3.10}$$

**Proof:**
(*a*) Analogous to the derivation of (2.6), it was shown in Theorem 5.2, Hunter, (1982) that *if* $G$ is any g-inverse of $I - P$, then for an irreducible *MRP*, the mean first passage time matrix is given by

$$M = \left[(1/\lambda)\{GP^{(1)}\Pi - E(GP^{(1)}\Pi)_d\} + I - G + EG_d\right]D, \tag{3.11}$$

where $\Pi = \boldsymbol{e\pi}^T$, $D = \lambda(\Pi_d)^{-1}$. Now $GP^{(1)}\Pi = GP^{(1)}\boldsymbol{e\pi}^T = G\boldsymbol{\mu\pi}^T$. From (3.11) and (3.6), and, since for any square matrix $X$, $tr(X) = \boldsymbol{e}^T X_d \boldsymbol{e}$, (3.6) follows.

(*b*) Taking $G = Z$, and since $Z\boldsymbol{e} = \boldsymbol{e}$, (3.7) leads to (3.8).

(c) Taking $G$ as the group inverse $A^{\#} = Z - \Pi$, and since $A^{\#}\boldsymbol{e} = \boldsymbol{0}$, (3.7) leads to (3.9).
Equations (3.8) and (3.9) are equivalent since $A^{\#}\boldsymbol{\mu} = Z\boldsymbol{\mu} - \lambda\boldsymbol{e}$, $A^{\#}\boldsymbol{\mu\pi}^T = Z\boldsymbol{\mu\pi}^T - \lambda\Pi$, and $tr(A^{\#}\boldsymbol{\mu\pi}^T)\boldsymbol{e} = tr(Z\boldsymbol{\mu\pi}^T)\boldsymbol{e} - \lambda\boldsymbol{e}$ and $tr(A^{\#}) = tr(Z) - 1$.

(*d*) With $\widetilde{G} \equiv [I - P + \boldsymbol{\mu u}^T]^{-1}$ (with $\boldsymbol{u}^T\boldsymbol{e} \neq 0$), Corollary 5.2.2, Hunter (1982) shows that

$$M = \left[I - \widetilde{G} + E\widetilde{G}_d\right]D, \tag{3.12}$$

where $D = \left[\boldsymbol{e}\boldsymbol{u}^T\widetilde{G}_d\right]^{-1}$. With this choice of g-inverse, using either (3.7) or (3.12), expression (3.10) for Kemeny's function vector $\boldsymbol{k}^{(1)}$ follows.

[Alternatively, use the result of Eqn. (3.17) of Hunter (1983): $[I-P+\boldsymbol{t}\boldsymbol{u}^T]^{-1}\boldsymbol{t} = (1/\boldsymbol{u}^T\boldsymbol{e})\boldsymbol{e}$, (provided $\boldsymbol{u}^T\boldsymbol{e} \neq 0$). Thus $\widetilde{G}\boldsymbol{\mu} = [I-P+\boldsymbol{\mu}\boldsymbol{u}^T]^{-1}\boldsymbol{\mu} = f\boldsymbol{e}$ with $f = (1/\boldsymbol{u}^T\boldsymbol{e})$ and $\widetilde{G}\boldsymbol{\mu}\boldsymbol{\pi}^T = f\Pi$. This leads to the first two terms of (3.7) cancelling, resulting in (3.10).]

**Theorem 2:** For a finite irreducible *MRP* with underlying irreducible *SMP* process having limiting vector $\boldsymbol{\varpi}$, mean sojourn time vector $\boldsymbol{\mu}$ and $\lambda \equiv \boldsymbol{\pi}^T\boldsymbol{\mu}$, Kemeny's function vector, $\boldsymbol{k}^{(2)}$ under Definition 2, can be expressed as follows:

(*a*) Let $G$ be any g-inverse of $I-P$, then

$$\boldsymbol{k}^{(2)} = \left[I + EG_d - (1/\lambda)E(G\boldsymbol{\mu}\boldsymbol{\pi}^T)_d\right]\boldsymbol{\mu}. \tag{3.13}$$

(*b*) If $Z = [I-P+\Pi]^{-1}$ is the fundamental matrix of the embedded *MC*,

$$\boldsymbol{k}^{(2)} = \left[I + EZ_d - (1/\lambda)E(Z\boldsymbol{\mu}\boldsymbol{\pi}^T)_d\right]\boldsymbol{\mu}. \tag{3.14}$$

(*c*) If $A^{\#}$ is the group inverse of $I-P$,

$$\boldsymbol{k}^{(2)} = \left[I + EA^{\#}_d - (1/\lambda)E(A^{\#}\boldsymbol{\mu}\boldsymbol{\pi}^T)_d\right]\boldsymbol{\mu}. \tag{3.15}$$

(*d*) If $\widetilde{G} \equiv [I-P+\boldsymbol{\mu}\boldsymbol{u}^T]^{-1}$ (with $\boldsymbol{u}^T\boldsymbol{e} \neq 0$),

$$\boldsymbol{k}^{(2)} = \boldsymbol{\mu} - (1/\boldsymbol{u}^T\boldsymbol{e})\boldsymbol{e} + \left(\boldsymbol{e}^T\widetilde{G}_d\boldsymbol{\mu}\right)\boldsymbol{e}. \tag{3.16}$$

**Proof:**

(*a*) Note that $\boldsymbol{k}^{(2)} = M\boldsymbol{\varpi}$ and, from (3.3), on taking the transpose, $\boldsymbol{\varpi} = \Lambda\boldsymbol{\pi}/\lambda$. Thus

$$D\boldsymbol{\varpi} = D\Lambda\boldsymbol{\pi}/\lambda = \frac{1}{\lambda}\begin{pmatrix} m_{11} & & \\ & m_{ii} & \\ & & m_{mm}\end{pmatrix}\begin{pmatrix} \mu_1 & & \\ & \mu_i & \\ & & \mu_m\end{pmatrix}\begin{bmatrix}\pi_1\\ \pi_i \\ \pi_m\end{bmatrix} = \frac{1}{\lambda}\begin{pmatrix} m_{11}\mu_1\pi_1 \\ m_{ii}\mu_i\pi_i \\ m_{mm}\mu_m\pi_m\end{pmatrix} = \boldsymbol{\mu}, \tag{3.17}$$

since from (3.5), $\lambda = \pi_i m_{ii}$ for all $i$.

Now from (3.11) and the definition of $\lambda$,

$$\boldsymbol{k}^{(2)} = M\boldsymbol{\varpi} = \left[(1/\lambda)\{G\boldsymbol{\mu}\boldsymbol{\pi}^T - E(G\boldsymbol{\mu}\boldsymbol{\pi}^T)_d\} + I - G + EG_d\right]\boldsymbol{\mu},$$

leading to (3.13).

(*b*) and (*c*): Equations (3.14) and (3.15) follow from (3.13).

(*d*): From (3.11), $\boldsymbol{k}^{(2)} = M\boldsymbol{\varpi} = M = \left[I - \widetilde{G} + E\widetilde{G}_d\right]D\boldsymbol{\varpi} = \left[I - \widetilde{G} + E\widetilde{G}_d\right]\boldsymbol{\mu}$

$= \boldsymbol{\mu} - \widetilde{G}\boldsymbol{\mu} + E\widetilde{G}_d\boldsymbol{\mu}$ where $\widetilde{G}\boldsymbol{\mu} = (1/\boldsymbol{u}^T\boldsymbol{e})\boldsymbol{e}$, from the proof of Theorem 1, leading to (3.16).

**Theorem 3:** For a finite irreducible *MRP* with underlying irreducible jump *MC* having transition matrix $P$ with stationary probability vector $\boldsymbol{\pi}$, mean sojourn time vector $\boldsymbol{\mu}$ and $\lambda \equiv \boldsymbol{\pi}^T\boldsymbol{\mu}$, Kemeny's function vector, $\boldsymbol{k}^{(3)}$ under Definition 3, can be expressed as follows:

(*a*) Let $G$ be any g-inverse of $I-P$, then

$$\boldsymbol{k}^{(3)} = (1/\lambda)\{G\boldsymbol{\mu} - tr(G\boldsymbol{\mu}\boldsymbol{\pi}^T)\}\boldsymbol{e} + \boldsymbol{e} - G\boldsymbol{e} + tr(G)\boldsymbol{e}. \tag{3.18}$$

($b$) If $Z = \left[I - P + \Pi\right]^{-1}$ is the fundamental matrix of the embedded $MC$,

$$\boldsymbol{k}^{(3)} = (1/\lambda\ )\{Z\boldsymbol{\mu} - tr(Z\boldsymbol{\mu}\boldsymbol{\pi}^T)\boldsymbol{e}\} + tr(Z)\boldsymbol{e}. \qquad (3.19)$$

($c$) If $A^{\#}$ is the group inverse of I – $P$,

$$\boldsymbol{k}^{(3)} = (1/\lambda)\{A^{\#}\boldsymbol{\mu} - tr(A^{\#}\boldsymbol{\mu}\boldsymbol{\pi}^T)\boldsymbol{e}\} + \boldsymbol{e} + tr(A^{\#})\boldsymbol{e}. \qquad (3.20)$$

($d$) If $\widetilde{G} \equiv [I - P + \boldsymbol{\mu}\boldsymbol{u}^T]^{-1}$ (with $\boldsymbol{u}^T\boldsymbol{e} \neq 0$)

$$\boldsymbol{k}^{(3)} = \left[I - \widetilde{G} + tr(\widetilde{G})\right]\boldsymbol{e}. \qquad (3.21)$$

**Proof:** From (3.6) observe that $\left(M_d\right)^{-1} = (1/\lambda)\Pi_d$ so that

$$\boldsymbol{k}^{(3)} = M\left(M_d\right)^{-1}\boldsymbol{e} = (1/\lambda)M\boldsymbol{\pi} = (1/\lambda)\boldsymbol{k}^{(1)}. \qquad (3.22)$$

Results (a), (b) (c) and (d) now follow directly from Theorem 1.

We now explore conditions under which the respective Kemeny function vectors $\boldsymbol{k}^{(i)}$ are multiples of $\boldsymbol{e}$, implying the constant nature of the $k_i^{(l)}$ and hence the existence of Kemeny constants for $MRP$'s, as is the case for Kemeny's constant for $MC$'s.

First note that

$(I - P)\boldsymbol{k}^{(1)} = (I - P)M\boldsymbol{\pi} = (P^{(1)}E - PM_d)\boldsymbol{\pi} = P^{(1)}\boldsymbol{e}\boldsymbol{e}^T\boldsymbol{\pi} - PD\boldsymbol{\pi} = P^{(1)}\boldsymbol{e} - P\lambda\boldsymbol{e} = \boldsymbol{\mu} - \lambda\boldsymbol{e}.$

Now $k_i^{(1)} = K_C^{(1)} \Leftrightarrow \boldsymbol{k}^{(1)} = K_C^{(1)}\boldsymbol{e} \Leftrightarrow (I - P)\boldsymbol{k}^{(1)} = \boldsymbol{0} \Leftrightarrow \boldsymbol{\mu} = \lambda\boldsymbol{e}.$

This is equivalent to $\mu_i = \lambda$. i.e. $\mu_i$ is a constant. (Note from its definition $\lambda = \sum_{k=1}^{m}\pi_k\mu_k$.)

This leads to a key result of this note:

**Corollary 1.1**: In an irreducible $MRP$ on a finite state space $S$, the mixture $k_i^{(1)} = \sum_{j\in S}\pi_j m_{ij}$, of the mean first passage times $m_{ij} = E[T_{ij} \mid X_0 = i]$ with the stationary distribution $\{\pi_j\}$ is independent of the initial state $i$ if and only if the mean sojourn times $\mu_j$ of visits to the state $j$ are the same for all states $j$.

When these conditions are satisfied, i.e. $\boldsymbol{\mu} = \lambda\boldsymbol{e}$, $\boldsymbol{k}^{(1)} = K_C^{(1)}\boldsymbol{e}$ where, under the conditions and notation of Theorem 1,

$K_C^{(1)} = \lambda\left[1 + tr(G) - tr(G\Pi)\right] = \lambda\left[tr(Z)\right] = \lambda\left[1 + tr(A^{\#})\right] = \lambda\left[1 + tr(\widetilde{G})\right] - (1/\boldsymbol{u}^T\boldsymbol{e}).$

The last expression follows from the observation that with $\widetilde{G} \equiv [I - P + \lambda\boldsymbol{e}\boldsymbol{u}^T]^{-1}$ then $\widetilde{G}\boldsymbol{e} \equiv [I - P + \lambda\boldsymbol{e}\boldsymbol{u}^T]^{-1}\boldsymbol{e} = (1/\lambda\boldsymbol{u}^T\boldsymbol{e})\boldsymbol{e}$, (provided $\boldsymbol{u}^T\boldsymbol{e} \neq 0$).

All of these results lead to a very simple expression for Kemeny's constant $K_C^{(1)}$, with the constant as derived for the embedded $MC$ case multiplied by the constant mean sojourn time.

Let us now consider the implication of the second definition of Kemeny's constant, $\boldsymbol{k}^{(2)} = M\boldsymbol{\varpi} = M\Lambda\boldsymbol{\pi}/\lambda$.

From equation (3.4),

$(I - P)\boldsymbol{k}^{(2)} = P^{(1)}\boldsymbol{e}\boldsymbol{e}^T\boldsymbol{\varpi} - PM_d\boldsymbol{\varpi} = \boldsymbol{\mu} - PM_d\Lambda\boldsymbol{\pi}/\lambda,$

where $M_d \Lambda \boldsymbol{\pi} = \begin{pmatrix} m_{11} & & \\ & m_{ii} & \\ & & m_{mm} \end{pmatrix} \begin{pmatrix} \mu_1 & & \\ & \mu_i & \\ & & \mu_m \end{pmatrix} \begin{pmatrix} \pi_1 \\ \pi_i \\ \pi_m \end{pmatrix} = \begin{pmatrix} m_{11}\mu_1\pi_1 \\ m_{ii}\mu_i\pi_i \\ m_{mm}\mu_m\pi_m \end{pmatrix} = \begin{pmatrix} \lambda\mu_1 \\ \lambda\mu_i \\ \lambda\mu_m \end{pmatrix}$

$= \lambda\boldsymbol{\mu},$ since, from (3.5), $\lambda = \pi_i m_{ii}$ for all $i$.

Now $(I - P)\boldsymbol{k}^{(2)} = \boldsymbol{\mu} - P\boldsymbol{\mu} = (I - P)\boldsymbol{\mu}$ i.e. $(I - P)(\boldsymbol{k}^{(2)} - \boldsymbol{\mu}) = \boldsymbol{0}.$

From (2.2) $(I - P)\boldsymbol{x} = \boldsymbol{0} \Leftrightarrow \boldsymbol{x} = c\boldsymbol{e}$ for some real $c$, so that $(I - P)(\boldsymbol{k}^{(2)} - \boldsymbol{\mu}) = \boldsymbol{0}$ if and only if $\boldsymbol{k}^{(2)} - \boldsymbol{\mu} = c\boldsymbol{e}$ for some real $c$, i.e. $\boldsymbol{k}^{(2)} = \boldsymbol{\mu} + c\boldsymbol{e}$.

Now $\boldsymbol{k}^{(2)} = K_C^{(2)}\boldsymbol{e} \Leftrightarrow K_C^{(2)}\boldsymbol{e} = \boldsymbol{\mu} + c\boldsymbol{e}$ i.e. if and only if $\boldsymbol{\mu}$ is a multiple of $\boldsymbol{e}$. Since $\lambda = \sum_{k=1}^{m} \pi_k \mu_k$ and all the $\mu_i$ are constant we must have $\mu_i = \lambda$. This leads to the following key result analogous to Corollary 1.1.

**Corollary 2.1**. In an irreducible *MRP* on a finite state space *S*, the mixture $k_i^{(2)} = \sum_{j \in S} \varpi_j m_{ij}$, of the mean first passage times $m_{ij} = E[T_{ij} \mid X_0 = i]$ with the stationary distribution $\{\varpi_j\}$ of the associated *SMP* is independent of the initial state *i* if and only if the mean sojourn times $\mu_j$ of visits to the state *j* are the same for all states *j*.

When these conditions are satisfied, i.e. $\boldsymbol{\mu} = \lambda\boldsymbol{e}$, $\boldsymbol{k}^{(2)} = K_C^{(2)}\boldsymbol{e}$ where, under the conditions and notation of Theorem 2,

$$K_C^{(2)} = \lambda\left[1 + tr(G) - tr(G\Pi)\right] = \lambda\left[tr(Z)\right] = \lambda\left[1 + tr(A^{\#})\right] = \lambda\left[1 + tr(\widetilde{G})\right] - (1/\boldsymbol{u}^T\boldsymbol{e}) = K_C^{(1)}.$$

Now, from (3.20), $\boldsymbol{k}^{(3)} = (1/\lambda)\boldsymbol{k}^{(1)}$. Directly from Corollary 1.1 we can derive an equivalent key result for $\boldsymbol{k}^{(3)}$.

**Corollary 3.1**: In an irreducible *MRP* on a finite state space *S*, the mixture $k_i^{(3)} = \sum_{j=1}^{m} m_{ij}/m_{jj}$ where $m_{ij} = E[T_{ij} \mid X_0 = i]$ are the mean first passage times between states *i* and *j* in the *MRP*, is independent of the initial state *i* if and only if the mean sojourn times $\mu_j$ of visits to the state *j* are the same for all states *j*.

Under the restriction $\boldsymbol{\mu} = \lambda\boldsymbol{e}$, with the notation and conditions of Theorem 3,

$$K_C^{(3)} = 1 + tr(G) - tr(G\Pi) = tr(Z) = 1 + tr(A^{\#}) = 1 + tr(\widetilde{G}) - (1/\lambda\boldsymbol{u}^T\boldsymbol{e}).$$

Analogous to the variant of Kemeny's constant for discrete time Markov chains we define the variants of the three Kemeny functions, as in Definitions 1, 2 and 3, as $k_i^{\circ(1)} \equiv \sum_{j=1, j\neq i}^{m} \pi_j m_{ij}$, $k_i^{\circ(2)} \equiv \sum_{j=1, j\neq i}^{m} \varpi_j m_{ij}$, and $k_i^{\circ(3)} \equiv \sum_{j=1, j\neq i}^{m} m_{ij}/m_{jj}$, respectively.

In Bini et. al. (2017) the notation for the "first passage time" and the "first return time" are defined as $\theta_i = \inf\{t \geq 0 : X(t) = i\}$ and $T_i = \inf\{t \geq J_1 : X(t) = i\}$, where $J_1$ ($T_1$ in our formulation of the definition of a *MRP*, $\{X_n, T_n\}$) is the first jump time of the *MC* so that if $X_0 = i$ then $\theta_i = 0 < T_i$ otherwise $\theta_i = T_i > 0$. Thus for $i \neq j, E_i[\theta_j] = E_i[T_j] = m_{ij}$ while for $i = j$, $E_i[\theta_i] = 0$ with $E_i[T_i] = m_{ii}$. Thus the deletion of $i = j$ term in the above sums retains just the $m_{ij}$ $= E_i[\theta_j]$ terms.

Since $\pi_i m_{ii} = \lambda$, (from (3.5)), $\varpi_i m_{ii} = \mu_i$, (from (3.3) or (3.17)), and $m_{ii}/m_{ii} = 1$, we deduce expressions for the equivalent vector functions as $\boldsymbol{k}^{\circ(1)} = \boldsymbol{k}^{(1)} - \lambda \boldsymbol{e}$, $\boldsymbol{k}^{\circ(2)} = \boldsymbol{k}^{(2)} - \boldsymbol{\mu}$, and $\boldsymbol{k}^{\circ(3)} = \boldsymbol{k}^{(3)} - \boldsymbol{e}$.
An immediate consequence of these alternative variants is that we can immediately deduce a variety of expressions for $\boldsymbol{k}^{\circ(1)}$, $\boldsymbol{k}^{\circ(2)}$ and $\boldsymbol{k}^{\circ(3)}$ directly from Theorems 1, 2 and 3.

**Corollary 1.2:** For a finite irreducible *MRP* with underlying irreducible jump *MC* having transition matrix *P* with stationary probability vector $\boldsymbol{\pi}$, mean sojourn time vector $\boldsymbol{\mu}$ and $\lambda \equiv \boldsymbol{\pi}^T \boldsymbol{\mu}$, the alternative Kemeny's function vector, $\boldsymbol{k}^{\circ(1)}$ , can be expressed as follows:

(*a*) Let *G* be any g-inverse of *I* – *P,* then

$$\boldsymbol{k}^{\circ(1)} = G\boldsymbol{\mu} - tr(G\boldsymbol{\mu}\boldsymbol{\pi}^T)\boldsymbol{e} - \lambda G\boldsymbol{e} + \lambda tr(G)\boldsymbol{e}. \tag{3.23}$$

(*b*) If $Z = [I - P + \Pi]^{-1}$ is the fundamental matrix of the embedded *MC*,

$$\boldsymbol{k}^{\circ(1)} = Z\boldsymbol{\mu} - tr(Z\boldsymbol{\mu}\boldsymbol{\pi}^T)\boldsymbol{e} + \lambda(tr(Z) - 1)\boldsymbol{e}. \tag{3.24}$$

(*c*) If $A^{\#}$ is the group inverse of *I* – *P,*

$$\boldsymbol{k}^{\circ(1)} = A^{\#}\boldsymbol{\mu} - tr(A^{\#}\boldsymbol{\mu}\boldsymbol{\pi}^T)\boldsymbol{e} + \lambda tr(A^{\#})\boldsymbol{e}. \tag{3.25}$$

(*d*) If $\widetilde{G} \equiv [I - P + \boldsymbol{\mu}\boldsymbol{u}^T]^{-1}$ (with $\boldsymbol{u}^T\boldsymbol{e} \neq 0$)

$$\boldsymbol{k}^{\circ(1)} = \lambda\left[tr(\widetilde{G}) - \widetilde{G}\right]\boldsymbol{e}. \tag{3.26}$$

**Corollary 2.2:** For a finite irreducible *MRP* with underlying irreducible *SMP* process having limiting vector $\boldsymbol{\varpi}$, mean sojourn time vector $\boldsymbol{\mu}$ and $\lambda \equiv \boldsymbol{\pi}^T \boldsymbol{\mu}$, the alternative Kemeny's function vector, $\boldsymbol{k}^{\circ(2)}$ , can be expressed as $\boldsymbol{k}^{\circ(2)} = \boldsymbol{K}^{\circ(2)}_C \boldsymbol{e}$ , where

(*a*) Let *G* be any g-inverse of *I* – *P,* then

$$K^{\circ(2)}_C = \boldsymbol{e}^T\left(G_d - (1/\lambda)(G\boldsymbol{\mu}\pi^T)_d\right)\boldsymbol{\mu}. \tag{3.27}$$

(*b*) If $Z = [I - P + \Pi]^{-1}$ is the fundamental matrix of the embedded *MC*,

$$K^{\circ(2)}_C = \boldsymbol{e}^T\left(Z_d - (1/\lambda)(Z\boldsymbol{\mu}\pi^T)_d\right)\boldsymbol{\mu}. \tag{3.28}$$

(*c*) If $A^{\#}$ is the group inverse of *I* – *P,*

$$K^{\circ(2)}_C = \boldsymbol{e}^T\left(A^{\#}_d - (1/\lambda)(A^{\#}\boldsymbol{\mu}\pi^T)_d\right)\boldsymbol{\mu}. \tag{3.29}$$

(*d*) If $\widetilde{G} \equiv [I - P + \boldsymbol{\mu}\boldsymbol{u}^T]^{-1}$ (with $\boldsymbol{u}^T\boldsymbol{e} \neq 0$),

$$K^{\circ(2)}_C = (1/\boldsymbol{u}^T\boldsymbol{e}) + \left(\boldsymbol{e}^T\widetilde{G}_d\boldsymbol{\mu}\right). \tag{3.30}$$

Thus the sum $k^{\circ(2)}_i \equiv \sum_{j=1, j\neq i}^{m} \varpi_j m_{ij}$ is a constant $\boldsymbol{K}^{\circ(2)}_C$, with one of the equivalent forms given by (3.27), (3.28), (3.29), or (3.30), for all values *i*, analogous to the result $\boldsymbol{K}^{\circ}{}_C$ for discrete time *MC's*.

**Corollary 3.2:** For a finite irreducible *MRP* with underlying irreducible jump *MC* having transition matrix $P$ with stationary probability vector $\boldsymbol{\pi}$, mean sojourn time vector $\boldsymbol{\mu}$ and $\lambda \equiv \boldsymbol{\pi}^T \boldsymbol{\mu}$, the alternative Kemeny's function vector, $\boldsymbol{k}^{\circ(3)}$**,** can be expressed as follows:

($a$) Let $G$ be any g-inverse of $I - P$, then

$$\boldsymbol{k}^{\circ(3)} = (1/\lambda)\{G\boldsymbol{\mu} - tr(G\boldsymbol{\mu}\boldsymbol{\pi}^T)\}\boldsymbol{e} - G\boldsymbol{e} + tr(G)\boldsymbol{e}. \quad (3.31)$$

($b$) If $Z = [I - P + \Pi]^{-1}$ is the fundamental matrix of the embedded *MC*,

$$\boldsymbol{k}^{\circ(3)} = (1/\lambda\ )\{Z\boldsymbol{\mu} - tr(Z\boldsymbol{\mu}\boldsymbol{\pi}^T)\boldsymbol{e}\} + tr(Z)\boldsymbol{e} - \boldsymbol{e}. \quad (3.32)$$

($c$) If $A^{\#}$ is the group inverse of $I - P$,

$$\boldsymbol{k}^{\circ(3)} = \ (1/\lambda)\{A^{\#}\boldsymbol{\mu} - tr(A^{\#}\boldsymbol{\mu}\boldsymbol{\pi}^T)\boldsymbol{e}\} + tr(A^{\#})\boldsymbol{e}. \quad (3.33)$$

($d$) If $\widetilde{G} \equiv [I - P + \boldsymbol{\mu}\boldsymbol{u}^T]^{-1}$ (with $\boldsymbol{u}^T\boldsymbol{e} \neq 0$)

$$\boldsymbol{k}^{\circ(3)} = \left[tr(\widetilde{G}) - \widetilde{G}\right]\boldsymbol{e}. \quad (3.34)$$

The variant $k_i^{\circ(2)} \equiv \sum_{j=1, j\neq i}^{m} \varpi_j m_{ij}$ is the only Kemeny function that achieves a constant for all *MRP's,* which naturally includes *MC's* in continuous times, birth and death processes, results that were also established in Bini et.al. (2017).

The observations that $\boldsymbol{k}^{\circ(1)} = \boldsymbol{k}^{(1)} - \lambda\boldsymbol{e}$ and $\boldsymbol{k}^{\circ(3)} = \boldsymbol{k}^{(3)} - \boldsymbol{e}$ imply that the results of Corollaries 1.1 and 3.1 also hold for $\boldsymbol{k}^{\circ(1)}$ and $\boldsymbol{k}^{\circ(3)}$ when the mean sojourn (holding) times, $\mu_j$, of visits to the state $j$ are the same for all states $j$ with

$$K_C^{\circ(1)} = \lambda\left[tr(G) - tr(G\Pi)\right] = \lambda\left[tr(Z) - 1\right] = \lambda\left[tr(A^{\#})\right] = \lambda\left[tr(\widetilde{G})\right] - (1/\boldsymbol{u}^T\boldsymbol{e}).$$

and

$$K_C^{\circ(3)} = tr(G) - tr(G\Pi) = tr(Z) = tr(A^{\#}) = tr(\widetilde{G}) - (1/\lambda\boldsymbol{u}^T\boldsymbol{e}).$$

Let us now explore a direct way of deriving the properties of the first Kemeny's function, $\boldsymbol{k} \equiv \boldsymbol{k}^{(1)} = \left(\sum_{j\in S} \pi_j m_{ij}\right)$.

From equation (3.5) i.e. $\boldsymbol{k} = P\boldsymbol{k} + \boldsymbol{\mu} - \lambda\boldsymbol{e}$ by repeated substitution we see that

$$\boldsymbol{k} = (I + P + .. + P^{n-1})(\boldsymbol{\mu} - \lambda\boldsymbol{e}) + P^n\boldsymbol{k}, \quad (3.35)$$

which, since, for each $r = 0, 1, \ldots$ $P^r\boldsymbol{e} = \boldsymbol{e}$, $\lambda = \boldsymbol{\pi}^T\boldsymbol{\mu}$ so that $P^r(\boldsymbol{\mu} - \lambda\boldsymbol{e}) = P^r\boldsymbol{\mu} - \boldsymbol{e}\boldsymbol{\pi}^T\boldsymbol{\mu} = (P^r - \Pi)\boldsymbol{\mu}$ and (3.21) reduces to

$$\boldsymbol{k} = (I - \Pi)\boldsymbol{\mu} + \sum_{k=1}^{n-1}(P^k - \Pi)\boldsymbol{\mu} + P^n\boldsymbol{k}. \quad (3.36)$$

Under the conditions of irreducibility and aperiodicity $P^n \to \boldsymbol{e}\boldsymbol{\pi}^T = \Pi$ and (Hunter (1983)) $Z = I + \sum_{k=1}^{\infty}(P - \Pi)^k = I + \sum_{k=1}^{\infty}(P^k - \Pi) = [I - P + \Pi]^{-1}$ so taking the limit as $n \to \infty$ leads to

$$\boldsymbol{k} = (I - \Pi)\boldsymbol{\mu} + (Z - I)\boldsymbol{\mu} + \Pi\boldsymbol{k} = (Z - \Pi)\boldsymbol{\mu} + \Pi\boldsymbol{k} = Z\boldsymbol{\mu} - \lambda\boldsymbol{e} + \Pi\boldsymbol{k}. \quad (3.37)$$

While (3.37) gives a relationship satisfied by $\boldsymbol{k}$ it does not lead to an explicit expression for $\boldsymbol{k}$. It does however illustrate how the $Z$ matrix emerges as a key tool for finding expressions for $\boldsymbol{k}$. It also is related to the "deviation matrix" $D=\sum_{n\geq 0}(P^n - \boldsymbol{e}\boldsymbol{\pi}^T) = Z - \Pi = A^{\#}$ which is used in the paper by Bini et al. (2017) in the special case of $MC$s in discrete time, where $k_i \equiv \sum_{j\varepsilon S} m_{ij}\pi_j = K_C = \sum_j D_{jj}$ .

A $MC$ in discrete time is a special $MRP$ with constant unit holding times ($T_{n+1} - T_n = 1$) between transitions, yielding $\mu_i = 1$ for all $i$ and implying $\lambda = 1$ and $K_C^{(1)} = K_C^{(2)} = K_C^{(3)} = tr(Z)$, $K^{\circ(1)}_C = K^{\circ(2)}_C = K^{\circ(3)}_C = tr(Z) - 1$, a constant not depending on $i$.

For a general $MRP$ it is still possible for Kemeny's function to be a constant with sojourn times in the states being random variables as long as their mean sojourn times are constant.

**Example 1:** *Two-state MRP.*

We consider the special two-state $MRP$ with $S = \{1, 2\}$, the embedded MC $\{X_n\}$ with transition matrix $P = \begin{bmatrix} p_{11} & p_{12} \\ p_{21} & p_{22} \end{bmatrix}$, with mean sojourn time vector $\boldsymbol{\mu} = (\mu_1, \mu_2)$. The stationary probability vector $\boldsymbol{\pi}^T = (\pi_1, \pi_2) = \left(p_{21}/(p_{12}+p_{21}), p_{12}/(p_{12}+p_{21})\right)$.
The constant $\lambda \equiv \boldsymbol{\pi}^T\boldsymbol{\mu} = \pi_1\mu_1 + \pi_2\mu_2 = (\mu_1 p_{21} + \mu_2 p_{12})/(p_{12}+p_{21})$. This implies that $\lambda(p_{12}+p_{21}) = \mu_1 p_{21} + \mu_2 p_{12}$.

Observe that $$Z = [I - P + \boldsymbol{e}\boldsymbol{\pi}^T]^{-1} = \frac{1}{p_{12}+p_{21}}\begin{bmatrix} p_{21} + \dfrac{p_{12}}{p_{12}+p_{21}} & p_{12} - \dfrac{p_{12}}{p_{12}+p_{21}} \\ p_{21} - \dfrac{p_{21}}{p_{12}+p_{21}} & p_{12} + \dfrac{p_{21}}{p_{12}+p_{21}} \end{bmatrix},$$

with $tr(Z) = 1 + \dfrac{1}{(p_{12}+p_{21})} = \dfrac{(1+p_{21}+p_{12})}{(p_{12}+p_{21})}$. (See Example 7.3.4, Hunter (1983)).

The mean first passage time matrix, $M$ is given by
$$M = \begin{bmatrix} m_{11} & m_{12} \\ m_{21} & m_{22} \end{bmatrix} = \begin{bmatrix} \mu_1 + (\mu_2 p_{12}/p_{21}) & \mu_1/p_{12} \\ \mu_2/p_{21} & (\mu_1 p_{21}/p_{12}) + \mu_2 \end{bmatrix}.$$ (See Eqn. (49) of Hunter (2016)).

Thus Kemeny's function, under the first definition, is given as
$$\boldsymbol{k}^{(1)} = \begin{bmatrix} k_1^{(1)} \\ k_2^{(1)} \end{bmatrix} = M\boldsymbol{\pi} = \frac{1}{(p_{12}+p_{21})}\begin{bmatrix} \mu_1(1+p_{21}) + \mu_2 p_{12} \\ \mu_1 p_{21} + \mu_2(1+p_{12}) \end{bmatrix}.$$
Observe that $\boldsymbol{k}^{(1)} = K_C^{(1)}\boldsymbol{e}$ if and only if $\mu_1(1+p_{21}) + \mu_2 p_{12} = \mu_1 p_{21} + \mu_2(1+p_{12})$,
i.e. $\mu_1 = \mu_2\ (= \lambda)$, with $K_C^{(1)} = \lambda(1+p_{21}+p_{12})/(p_{12}+p_{21}) = \lambda tr(Z)$, a constant, as expected from Corollary 1.1 regarding the conditions for the existence of Kemeny's constant.
$$\boldsymbol{k}^{\circ(1)} = \begin{bmatrix} k^{\circ(1)}_1 \\ k^{\circ(1)}_2 \end{bmatrix} = \boldsymbol{k}^{(1)} - \lambda\boldsymbol{e} = \frac{1}{(p_{12}+p_{21})}\begin{bmatrix} \mu_1 \\ \mu_2 \end{bmatrix} = \begin{bmatrix} \pi_2 m_{12} \\ \pi_1 m_{21} \end{bmatrix}.$$

Further, the limiting probability vector associated with the *SMP* is $\boldsymbol{\varpi}^T = \boldsymbol{\pi}^T \Lambda/\lambda = \{\varpi_1, \varpi_2\}$ where, with $\Lambda = diag(\mu_1, \mu_2)$ yields $\boldsymbol{\varpi}^T = (\mu_1 p_{21}/(\mu_1 p_{21} + \mu_2 p_{12}), \mu_2 p_{12}/(\mu_1 p_{21} + \mu_2 p_{12}))$. Thus Kemeny's function, under the second definition,

$$\boldsymbol{k}^{(2)} = \begin{bmatrix} k_1^{(2)} \\ k_2^{(2)} \end{bmatrix} = M\boldsymbol{\varpi} = \frac{1}{\mu_1 p_{21} + \mu_2 p_{12}} \begin{bmatrix} \mu_1\{\mu_1 p_{21} + \mu_2(1+p_{12})\} \\ \mu_2\{\mu_1(1+p_{21}) + \mu_2 p_{12}\} \end{bmatrix}.$$

Observe that $\boldsymbol{k}^{(2)} = K_C^{(2)}\boldsymbol{e}$ if and only if $\mu_1 = \mu_2\ (=\lambda)$, $K_C^{(2)} = \lambda + \lambda/(p_{21}+p_{12}) = \lambda(1+p_{12}+p_{21})/(p_{12}+p_{21}) = \lambda tr(Z)$, a constant, as expected from Corollary 2.1

Note also that $\boldsymbol{k}^{\circ(2)} = \boldsymbol{k}^{(2)} - \boldsymbol{\mu} = \dfrac{\mu_1\mu_2}{\mu_1 p_{21} + \mu_2 p_{12}} \begin{bmatrix} 1 \\ 1 \end{bmatrix} = \begin{bmatrix} m_{12}\varpi_2 \\ m_{21}\varpi_1 \end{bmatrix}$, leading to a constant for each state, as expected , as a consequence of Corollary 2.2.

Further

$$\boldsymbol{k}^{(3)} = \begin{bmatrix} k_1^{(3)} \\ k_2^{(3)} \end{bmatrix} = (1/\lambda)\boldsymbol{k}^{(1)} = \frac{1}{(\mu_1 p_{21} + \mu_2 p_{12})} \begin{bmatrix} \mu_1(1+p_{21}) + \mu_2 p_{12} \\ \mu_1 p_{21} + \mu_2(1+p_{12}) \end{bmatrix}.$$

Observe that $\boldsymbol{k}^{(3)} = K_C^{(3)}\boldsymbol{e}$ if and only if $\mu_1 = \mu_2\ (=\lambda)$, with $K_C^{(3)} = (1+p_{21}+p_{12})/(p_{21}+p_{12})$ a constant, $tr(Z)$, as expected from Corollary 3.1 regarding the conditions for the existence of Kemeny's constant.
Further,

$$\boldsymbol{k}^{\circ(3)} = \begin{bmatrix} k^{\circ(3)}_1 \\ k^{\circ(3)}_2 \end{bmatrix} = \boldsymbol{k}^{(3)} - \boldsymbol{e} = \frac{1}{(\mu_1 p_{21} + \mu_2 p_{12})} \begin{bmatrix} \mu_1 \\ \mu_2 \end{bmatrix} = \begin{bmatrix} m_{12}/m_{22} \\ m_{21}/m_{11} \end{bmatrix}.$$

## 4. Kemeny's function for Markov chains in continuous time on a finite state space.

*MC's* in continuous time are special *MRP's* whose semi-Markov kernel is given by $Q_{ij}(t) = p_{ij}\left[1 - e^{-\nu_i t}\right]$, $(i, j = 1,..,m,\ t \ge 0)$ with $p_{ii} = 0$. We shall assume that $0 < \nu_i < \infty$ so that the process is stable and regular. In the terminology of Section 3, $F_{ij}(t) = 1 - e^{-\nu_i t}$, $i \ne j$, $t \ge 0$, implying that for $(i,j) \in S = \{1,2,...,m\}$, $\mu_{ij} = p_{ij}/\nu_i$ $(i \ne j)$ and $\mu_i = 1/\nu_i$, $i \in S$.

It is typical to specify a *MC* in continuous time in terms of the infinitesimal generator $Q = [q_{ij}]$ of the process. However, (see Çinlar (1975) or Hunter (1969)), $Q$ can be expressed in terms of the stated parameters from the *MRP*. In particular

$$q_{ij} = \begin{cases} -\nu_i, & i = j, \\ \nu_i p_{ij}, & i \ne j. \end{cases} \quad \text{Conversely, } \nu_i = -q_{ii}, \text{ and } p_{ij} = \begin{cases} 0, & i = j, \\ -q_{ij}/q_{ii}, & i \ne j. \end{cases}$$

Observe that $q_{ij} \geq 0$ for $i \neq j$, that $q_{ii} < 0$, and that $\sum_j q_{ij} = 0$ for all $i \in S$. It is easy to see that $I - P = (Q_d)^{-1} Q$ and $\Lambda = diag(\mu_1, \mu_2, ..., \mu_m) = -(Q_d)^{-1}$.

The *MRP* expressions for the Kemeny's functions, under the above definitions, also apply to *MC's* in continuous time. In fact, a *MC* in continuous time can also have a constant Kemeny's function if the mean sojourn times in each state are exponentially distributed with the same mean ν. This follows directly from the results of Section 3. If $\nu_i = \nu$, for all $i \in S$, then a Poisson process of rate ν is driving the *MC* in continuous time process with transition rates $q_{ij} = \nu p_{ij}$, $(i \neq j)$, or transition probabilities $p_{ij} = q_{ij} \big/ \sum_{k \neq i} q_{ik}$, $(i \neq j)$.

Let $\{X_t, t \geq 0\}$ be the associated *SMP* i.e. $X_t = X_n$ for $T_n \leq t < T_{n+1}$, the state of the *MRP* at time *t*. It is typical for regular *MC's* in continuous time to consider the limiting probabilities $\varpi_j = \lim_{t \to \infty} P\{X_t = j | X_0 = i\}$ rather than the stationary probabilities $\pi_i$ of the embedded *MC*.

Both of these sets of probabilities can be found as the solutions of stationary equations. Let $\boldsymbol{\varpi}^T = \{\varpi_1, \varpi_2, ..., \varpi_m\}$. It is easily seen that $\boldsymbol{\pi}^T(I - P) = 0 \Leftrightarrow \boldsymbol{\varpi}^T Q = 0^T$. Further, since $\lambda = \boldsymbol{\pi}^T \boldsymbol{\mu}$, $\boldsymbol{\varpi}^T = \boldsymbol{\pi}^T \Lambda / \lambda$ and $\boldsymbol{\pi}^T = \lambda \boldsymbol{\varpi}^T \Lambda^{-1}$.

Through the identification of the transition probabilities $p_{ij}$, the mean holding times $\mu_i$, and a generalized inverse *G* of *I* – *P* where *P* = $[p_{ij}]$ one can use one of the results of Theorem 4 to find an expression for Kemeny's function, under Definition 1. However it is typical to prescribe the *MC* in continuous time simply in terms of the infinitesimal generator *Q* = $[q_{ij}]$ so it is desirable to express Kemeny's function in terms of *Q*.

Since Kemeny's function, under Definition 1, is $\boldsymbol{k}^{(1)} = M\boldsymbol{\pi}$ we can use the results of Theorem 5.4 of Hunter (1982) where an expression for *M* is derived in terms of $H = \left[Q + \boldsymbol{e}\boldsymbol{u}^T\right]^{-1}$ (with $\boldsymbol{u}^T\boldsymbol{e} \neq 0$), a generalized inverse of *Q*.
We do not go into the details but refer the reader to Hunter (1982) where it is shown that
$\boldsymbol{\varpi}^T = \boldsymbol{u}^T H$, $\lambda = -\left(\boldsymbol{\varpi}^T Q_d \boldsymbol{e}\right)^{-1}$, $\boldsymbol{\pi}^T = -\lambda \boldsymbol{\varpi}^T Q_d$, and $M = \left[H - EH_d - (Q_d)^{-1}\right]\left[(\boldsymbol{e}\boldsymbol{u}^T H)_d\right]^{-1}$.
Further $D = M_d = \lambda (\Pi_d)^{-1}$ where $\Pi_d = -\lambda \left(\boldsymbol{e}\boldsymbol{u}^T H\right)_d Q_d$ leading to $D = -(Q_d)^{-1}\left[\left(\boldsymbol{e}\boldsymbol{u}^T H\right)_d\right]^{-1}$.
Since $D\boldsymbol{\pi} = \lambda \boldsymbol{e}$, $\boldsymbol{k}^{(1)} = M\boldsymbol{\pi} = \lambda\left[I - HQ_d + EH_d Q_d\right]$.
i.e.

$$\boldsymbol{k}^{(1)} = \lambda\left[I - HQ_d + tr(HQ_d)\right]\boldsymbol{e}. \tag{4.1}$$

This last result also follows from (3.12) since $\widetilde{G} = [I - P + (Q_d)^{-1}\boldsymbol{e}\boldsymbol{u}^T]^{-1} = HQ_d$.

This result (4.1) maybe not be particularly useful, as it does require the computation of a matrix inverse. However it reflects the similarity between the computations for $\boldsymbol{k}$ by using transition rates rather than the transition probabilities.

**Example 2:** *Two-state MC in continuous time*

For the two-state *MC* in continuous time with $Q = \begin{bmatrix} -\nu_1 & \nu_1 \\ \nu_2 & -\nu_2 \end{bmatrix}$, we have $P = \begin{bmatrix} 0 & 1 \\ 1 & 0 \end{bmatrix}$ so that the underlying *MRP* $\{(X_n, T_n), n \geq 1\}$ process is effectively an alternating renewal process $\{U_1, V_1, U_2, V_2, \ldots\}$ with $U_i = T_{2i-1} - T_{2i-2}$, $V_i = T_{2i} - T_{2i-1}$, $(i \geq 1)$ where $T_0 = 0$. Thus $U_i$ is distributed as an exponential ($\nu_1$) random variable and $V_i$ is distributed as an exponential ($\nu_2$) random variable implying that $\mu_1 = 1/\nu_1$, $\mu_2 = 1/\nu_2$.

Since $\boldsymbol{\pi} = \begin{pmatrix} \pi_1 \\ \pi_2 \end{pmatrix} = \begin{pmatrix} 1/2 \\ 1/2 \end{pmatrix}$, $M = \begin{bmatrix} m_{11} & m_{12} \\ m_{21} & m_{22} \end{bmatrix} = \begin{bmatrix} 1/\nu_1 + 1/\nu_2 & 1/\nu_1 \\ 1/\nu_2 & 1/\nu_1 + 1/\nu_2 \end{bmatrix}$,
it follows that Kemeny's function, under Definition 1, is given by

$$\boldsymbol{k}^{(1)} = M\boldsymbol{\pi} = \begin{pmatrix} (1/\nu_1) + (1/2\nu_2) \\ (1/2\nu_1) + (1/\nu_2) \end{pmatrix}.$$

Thus $\boldsymbol{k}^{(1)} = k\boldsymbol{e} \Leftrightarrow \nu_1 = \nu_2 = \nu$, say, with $k = 3/2\nu$, in which case the two-state $MC$ in continuous time reduces to a Renewal process with exponential holding times, i.e. a Poisson process.

Since $\boldsymbol{\varpi} = \begin{pmatrix} \varpi_1 \\ \varpi_2 \end{pmatrix} = \begin{pmatrix} \nu_2/(\nu_1+\nu_2) \\ \nu_1/(\nu_1+\nu_2) \end{pmatrix}$, it follows that Kemeny's function under Definition 2, is given by $\boldsymbol{k}^{(2)} = M\boldsymbol{\varpi} = \begin{pmatrix} (1/\nu_1) + (1/(\nu_1+\nu_2)) \\ (1/\nu_2) + (1/(\nu_1+\nu_2)) \end{pmatrix}$

Similarly $\boldsymbol{k}^{(2)} = k\boldsymbol{e} \Leftrightarrow \nu_1 = \nu_2 = \nu$, say, with $k = 3/2\nu$, as above.

$\boldsymbol{k}^{o(2)} = M\boldsymbol{\varpi} - \boldsymbol{\mu} = \dfrac{1}{\nu_1+\nu_2}\begin{pmatrix} 1 \\ 1 \end{pmatrix} = \begin{bmatrix} m_{12}\varpi_2 \\ m_{21}\varpi_1 \end{bmatrix}$, as expected from Corollary 2.2.

In conclusion, while under special situations Markov renewal processes including discrete Markov chains in continuous time, may retain the special constant feature of Kemeny's constant exhibited in discrete time *MC's*, the constant nature does not hold in general.

We conclude by considering a general three-state Birth and Death process.

**Example 3:** *Three-state Birth and Death process in continuous time*

We consider the process with state-space S = {1, 2, 3}, with birth rates $\alpha_i$, for $i$ = 1, 2 and death rates $\beta_i$ = 2, 3. Let $\rho_2 = \alpha_2/\beta_2$.

Thus $Q = \begin{bmatrix} -\alpha_1 & \alpha_1 & 0 \\ \beta_2 & -(\alpha_2+\beta_2) & \alpha_2 \\ 0 & \beta_3 & -\beta_3 \end{bmatrix}$, $P = \begin{bmatrix} 0 & 1 & 0 \\ 1/(1+\rho_2) & 0 & \rho_2/(1+\rho_2) \\ 0 & 1 & 0 \end{bmatrix}$, are respectively the infinitesimal generator of the $MC$ in continuous time and the transition matrix of the embedded discrete time $MC$.

The stationary probability vector of the embedded $MC$ is given by the solution of $\boldsymbol{\pi}^T = \boldsymbol{\pi}^T P$ giving $\boldsymbol{\pi}^T = (\pi_1, \pi_2, \pi_3) = \left( 1/2(1+\rho_2),\ 1/2,\ \rho_2/2(1+\rho_2) \right)$.
The stationary probability vector of the continuous-time $MC$ is given by the solution of $\boldsymbol{\varpi}^T Q = \boldsymbol{0}^T$, giving $\boldsymbol{\varpi}^T = (\varpi_1, \varpi_2, \varpi_3) = \left( \varpi_1,\ (\alpha_1/\beta_2)\varpi_1\ \ (\alpha_1\alpha_2/\beta_2\beta_3)\varpi_1 \right)$ where $\varpi_1 = 1/(1 + \alpha_1/\beta_2 + \alpha_1\alpha_2/\beta_2\beta_3)$.

Observe that $P^{(1)} = \left[\mu_{ij}\right] = \begin{bmatrix} 0 & 1/\alpha_1 & 0 \\ (1/1+\rho_2)(1/\alpha_2+\beta_2) & 0 & (\rho_2/1+\rho_2)(1/\alpha_2+\beta_2) \\ 0 & 1/\beta_3 & 0 \end{bmatrix}$, leading to the vector of mean holding times $\boldsymbol{\mu} = \begin{bmatrix} \mu_1 \\ \mu_2 \\ \mu_3 \end{bmatrix} = \begin{bmatrix} 1/\alpha_1 \\ 1/(\alpha_2+\beta_2) \\ 1/\beta_3 \end{bmatrix}$.

The mean asymptotic increment. $\lambda = \boldsymbol{\pi}^T \boldsymbol{\mu} = \dfrac{1}{2(1+\rho_2)}\left[\dfrac{1}{\alpha_1} + \dfrac{1}{\beta_2} + \dfrac{\alpha_2}{\beta_2\beta_3}\right] = \dfrac{1}{2(1+\rho_2)\alpha_1\varpi_1}$.

In this special case one can solve for the mean first passage times directly from equation (3.2) or (3.4). This leads to following expressions

(1) $m_{11} = \dfrac{1}{\alpha_1} + \dfrac{1}{\beta_2} + \dfrac{\rho_2}{\beta_3} = 2\lambda(1+\rho_2)$

(2) $m_{12} = \dfrac{1}{\alpha_1}$,

(3) $m_{13} = \dfrac{1}{\alpha_1} + \dfrac{1}{\alpha_2} + \dfrac{1}{\alpha_1\rho_2} = m_{33} + \dfrac{1}{\alpha_1}$,

(4) $m_{21} = \dfrac{1}{\beta_2} + \dfrac{\rho_2}{\beta_3} = m_{11} - \dfrac{1}{\alpha_1}$,

(5) $m_{22} = \dfrac{1}{(1+\rho_2)}\left(\dfrac{1}{\alpha_1} + \dfrac{1}{\beta_2} + \dfrac{\rho_2}{\beta_3}\right) = \dfrac{m_{11}}{(1+\rho_2)} = 2\lambda$,

(6) $m_{23} = \dfrac{1}{\rho_2\alpha_1} + \dfrac{1}{\alpha_2} = m_{33} - \dfrac{1}{\beta_2}$

(7) $m_{31} = \dfrac{1}{\beta_2} + \dfrac{1+\rho_2}{\beta_3} = m_{11} - \dfrac{1}{\alpha_1} + \dfrac{1}{\beta_3}$,

(8) $m_{32} = \dfrac{1}{\beta_3}$

(9) $m_{33} = \dfrac{1}{\rho_2\alpha_1} + \dfrac{1}{\alpha_2} + \dfrac{1}{\beta_3} = \dfrac{(1+\rho_2)m_{22}}{\rho_2} = \dfrac{m_{11}}{\rho_2} = \dfrac{2\lambda(1+\rho_2)}{\rho_2}$

Note that simplification of expressions can be effected by noting that $2\lambda(1+\rho_2)\alpha_1\varpi_1 = 1$.
Since $m_{11} = 2\lambda(1+\rho_2)$, $m_{22} = 2\lambda$, $m_{33} = 2\lambda(1+\rho_2)/\rho_2$ this implies that
$m_{11}\varpi_1 = 1/\alpha_1$, $m_{22}\varpi_2 = 1/(\alpha_2+\beta_2)$, $m_{33}\varpi_3 = 1/\beta_3$, or as a result that $m_{ii}\varpi_i = \mu_i$.

With these expressions we can now find expressions for the three Kemeny vector functions and their variants:

$$\boldsymbol{k}^{(1)} = \boldsymbol{M\pi} = \begin{bmatrix} k_1^{(1)} \\ k_2^{(1)} \\ k_3^{(1)} \end{bmatrix} = \begin{bmatrix} 2\lambda + \dfrac{\pi_1 + \pi_3}{\alpha_1} - \dfrac{\pi_3}{\beta_3} \\ 3\lambda - \dfrac{\pi_1}{\alpha_1} - \dfrac{\pi_3}{\beta_3} \\ 2\lambda - \dfrac{\pi_1}{\alpha_1} + \dfrac{(\pi_1 + \pi_2)}{\beta_3} \end{bmatrix}, \boldsymbol{k}^{\circ(1)} = \boldsymbol{M\pi} - \lambda \boldsymbol{e} = \begin{bmatrix} \lambda + \dfrac{\pi_1 + \pi_3}{\alpha_1} - \dfrac{\pi_3}{\beta_3} \\ 2\lambda - \dfrac{\pi_1}{\alpha_1} - \dfrac{\pi_3}{\beta_3} \\ \lambda - \dfrac{\pi_1}{\alpha_1} + \dfrac{(\pi_1 + \pi_2)}{\beta_3} \end{bmatrix}.$$

leading immediately to $k_1^{(1)} = k_2^{(1)} \Leftrightarrow \lambda = \dfrac{1}{\alpha_1}$; $k_2^{(1)} = k_3^{(1)} \Leftrightarrow \alpha_1 = \beta_3$; $k_1^{(1)} = k_3^{(1)} \Leftrightarrow \alpha_1 = \beta_3$;

From this it is easy to deduce that equality of the first type Kemeny functions holds when $\alpha_1 = \alpha_2 + \beta_2 = \beta_3$ consistent with the equality of the mean holding times $\mu_1 = \mu_2 = \mu_3$.

Similarly

$$\boldsymbol{k}^{(2)} = \boldsymbol{M\varpi} = \begin{bmatrix} k_1^{(2)} \\ k_2^{(2)} \\ k_3^{(2)} \end{bmatrix} = \begin{bmatrix} \dfrac{1}{\alpha_1} + \dfrac{\varpi_2}{\alpha_1} + \left(\dfrac{1}{\alpha_1} - \dfrac{1}{\beta_3}\right)\varpi_3 + \dfrac{1}{\beta_3} \\ \left(\dfrac{1}{\alpha_1} - \dfrac{\varpi_1}{\alpha_1}\right)_1 + \dfrac{1}{\alpha_2 + \beta_2} + \left(\dfrac{1}{\beta_3} - \dfrac{\varpi_3}{\beta_3}\right) \\ \left(\dfrac{1}{\alpha_1} - \dfrac{\varpi_1}{\alpha_1} + \dfrac{\varpi_1}{\beta_3}\right) + \dfrac{\varpi_2}{\beta_3} + \dfrac{1}{\beta_{3\,3}} \end{bmatrix}$$

implying $k_1^{(2)} = k_3^{(2)} \Leftrightarrow \alpha_1 = \alpha_2 + \beta_2$; $k_1^{(2)} = k_3^{(2)} \Leftrightarrow \alpha_1 = \beta_3$; $k_2^{(2)} = k_3^{(2)} \Leftrightarrow \alpha_2 + \beta_2 = \beta_3$.

Thus the equality of the second type of Kemeny functions also hold if and only if $\alpha_1 = \alpha_2 + \beta_2 = \beta_3$ as for the first type.

$$\boldsymbol{k}^{\circ(2)} = \boldsymbol{M\varpi} - \boldsymbol{\mu} = \begin{bmatrix} \dfrac{\varpi_2}{\alpha_1} + \left(\dfrac{1}{\alpha_1} - \dfrac{1}{\beta_3}\right)\varpi_3 + \dfrac{1}{\beta_3} \\ \left(\dfrac{1}{\alpha_1} - \dfrac{\varpi_1}{\alpha_1}\right)_1 + \left(\dfrac{1}{\beta_3} - \dfrac{\varpi_3}{\beta_3}\right) \\ \left(\dfrac{1}{\alpha_1} - \dfrac{\varpi_1}{\alpha_1} + \dfrac{\varpi_1}{\beta_3}\right) + \dfrac{\varpi_2}{\beta_3} \end{bmatrix} = \left(\dfrac{1 - \varpi_1}{\alpha_1} + \dfrac{1 - \varpi_3}{\beta_3}\right)\begin{bmatrix} 1 \\ 1 \\ 1 \end{bmatrix},$$

leading to a constant value for each term $k_i^{\circ(2)} = \sum_{j=1, j\neq i}^{3} \varpi_j m_{ij}$, as expected from Corollary2.2.

The equality conditions for the third type are as for the first type by virtue of the earlier observation that $\boldsymbol{k}^{(3)} = (1/\lambda)\boldsymbol{k}^{(1)}$.

In conclusion we have established and illustrated by a number of examples that for a finite irreducible Markov renewal process and the special cases of Markov chains in continuous time (including Birth – death processes) Kemeny's functions – defined in a variety ways - are constant if and only if the mean sojourn time of the *MRP* on any visit to a state is constant. However, for one particular variant, $\sum_{j=1, j\neq i}^{m} \varpi_j m_{ij}$, this sum is constant for all states in any finite irreducible *MRP*. This is consistent with the observations of Bini et. al. (2017) leading to a natural extension of Kemeny's constant to not only Markov chains in continuous time but to Markov renewal processes in general.

In this paper we have not examined the case of an infinite countable state space having used matrix techniques involving finite dimensional matrices to develop the results. The paper by Bini et. al, (2017) extends the results to countably infinite state spaces on not only discrete time Markov chains but also continuous time Markov chains. The arguments involve conditions for the convergence or divergence of Kemeny constant functions.

**Acknowledgement**
The author wishes to acknowledge the discussions and communications with Professor Daryl Daley, University of Melbourne in respect to the extension of Kemeny's constant to Markov renewal processes. These were initiated when they were both invited speakers at the conference on "Markov & Semi-Markov Processes & Related Fields" at Sithonia, Greece, September 20-23, 2011and have continued over the intervening period.
The extension to the variants of Kemeny's constants, as highlighted by Corollaries 1.2, 2.2 and 3.2, was a result of some helpful comments made by Guy Latouche.